\documentclass[leqno,english]{amsart}
\usepackage{amsmath,amssymb,amsthm,enumerate,esint,mathtools,xparse}
\usepackage[svgnames]{xcolor}
 \usepackage{placeins}
 \usepackage{float}
 \usepackage{subfig}
  \usepackage{hyperref}
 \usepackage{empheq}
\usepackage{cancel}

  \usepackage{bbm}

 \usepackage{mathrsfs} 
\usepackage{mathabx}

\numberwithin{equation}{section}

\newtheorem{theorem}{Theorem}[section]

\newtheorem{remark}[theorem]{Remark}

\newcommand{\bke}[1]{\left ( #1 \right )}
\newcommand{\bkt}[1]{\left [ #1 \right ]}

\newcommand{\norm}[1]{\left \| #1 \right \|}

\newcommand{\abs}[1]{\left | #1 \right |}

\newcommand\ga{\gamma}

\newcommand\ze{\zeta}
\renewcommand\th{\theta}
\newcommand\ka{\kappa}
\newcommand\la{\lambda}
\newcommand\si{\sigma}

\newcommand\De{\Delta}

\newcommand\Om{\Omega}
\newcommand\Rg{\mathscr{R}_g}

\newcommand{\R}{\mathbb{R}}

\renewcommand{\div}{\mathop{\rm div}\nolimits}
\newcommand{\curl} {\mathop{\rm curl}\nolimits}

\newcommand{\pd}{\partial}
\newcommand{\nb}{\nabla}

\newcommand{\EQ}[1]{\begin{equation}\begin{split} #1 \end{split}\end{equation}}
\newcommand{\EQN}[1]{\begin{equation*}\begin{split} #1 \end{split}\end{equation*}}

\DeclarePairedDelimiter{\oldnormaux}{\bracevert}{\bracevert}

\NewDocumentCommand{\oldnorm}{som}{%
  \IfBooleanTF{#1}
    {\oldnormaux*{#3}}
    {\IfNoValueTF{#2}
       {\oldnormaux*{\vphantom{dq}#3}}
       {\oldnormaux[#2]{#3}}%
    }%
}

\begin{document}
\title[Asymmetric equilibrium bubbles]{
Nonspherically symmetric equilibrium bubbles in a steadily rotating incompressible
fluid}

\author[C.-C. Lai]{Chen-Chih Lai}
\address{\noindent
Department of Mathematics,
Columbia University , New York, NY, 10027, USA}
\email{cclai.math@gmail.com}

\author[M. I. Weinstein]{Michael I. Weinstein}
\address{\noindent
Department of Applied Physics and Applied Mathematics and Department of Mathematics,
Columbia University , New York, NY, 10027, USA}
\email{miw2103@columbia.edu}

\begin{abstract}
This note presents two nontrivial, rotational equilibrium solutions to the spatial uniform gas pressure (isobaric) approximate model of Prosperetti in the inviscid case.
Building on Gavrilov's work [GAFA 2019], we first establish the existence of equilibrium solutions with nontrivial (rotational) liquid flow.
Second, we construct a nonspherically symmetric, horn-torus-shaped equilibrium bubble under mild spatial decay conditions of the liquid flow.
In addition, we extend earlier results on the characterization of spherical equilibrium bubbles to the axisymmetric, purely azimuthal setting.
Finally, we implement a numerical simulation of the equilibrium bubble shape using the Physics-Informed Neural Network (PINN) approximation.
\end{abstract}
\maketitle

\tableofcontents

\section{Introduction}

In this paper we study equilibrium solutions of the isobaric (uniform gas pressure) approximate model of Prosperetti  \cite{Prosperetti-JFM1991} for a gas bubble, in which thermal diffusion occurs, immersed  in an incompressible liquid with surface tension. We refer to this model throughout as ``the approximate model''. It  has been studied as well by  Biro--Vel\'azquez \cite{bv-SIMA2000}, and the authors \cite{LW-vbas2022, LW-vbaslinear2023, LW-ipug}.

Our main result is the construction of two nontrivial (nonspherically symmetric, rotational) equilibrium solutions to the approximate model.

In addition, we extend earlier characterization results on sphericity of equilibrium bubbles.
Prior results--namely, \cite[Proposition 4.3,(1)]{LW-vbas2022} and \cite[Theorem 3.1]{LW-ipug}--show that equilibrium bubbles must be spherical under the assumptions that either $\mu_l>0$ (nonzero viscosity in the liquid), or the equilibrium liquid velocity is irrotational, along with suitable far-field conditions.
In Part (a) of Theorem \ref{thm-main}, we extend these results to axisymmetric, purely azimuthal flows.

Finally, in Appendix \ref{sec:numerics}, we implement a numerical simulation of the equilibrium bubble shape by applying the Physics-Informed Neural Network (PINN) method \cite{PINN-JCP2019} to solve the stress balance equation (Laplace--Young condition).
The learned bubble surface is then compared with the analytical solution from Theorem \ref{thm-main}.
This demonstrates the potential of PINNs as a flexible tool for exploring equilibrium configurations beyond the analytically tractable regime.

\subsection{Background of the approximate model}

There are extensive studies of the approximate model in the spherically symmetric setting, e.g. the dynamics arising when a spherical equilibriium is subjected to a spherically symmetric perturbation.  
To investigate the thermal effects in the bubble oscillations, in the early 90s, Prosperetti \cite{Prosperetti-JFM1991} derived the approximate model of a spherically symmetric gas bubble deforming in an incompressible liquid assuming spatially uniform gas pressure.
In this paper, he also computed approximated values of the thermal dissipation rate for the nearly isothermal and the nearly adiabatic cases.
Almost a decade later, Biro and Vel\'azquez \cite{bv-SIMA2000} derived the similar model, motivated by the  regime of sonoluminescence experiments.
They analyzed the model mathematically and proved the local-in-time well-posedness of the spherically symmetric model using the Schauder fixed point theorem and the classical regularity theory for quasilinear parabolic equations.
They established the Lyapunov stability of the spherically symmetric equilibrium relative to small spherically symmetric perturbations.
In  \cite{LW-vbas2022}, the authors went beyond Lyapunov stability  to establish exponential asymptotic stability.
More precisely, we proved that the spherically symmetric model admits a manifold of equilibria, parametrized by the bubble mass, and that this manifold of equilibria is exponentially asymptotically stable relative to small spherically symmetric perturbations.
In \cite{LW-vbaslinear2023}, we examined the periodically forced problem, where the dynamics of the bubble is driven by a small-amplitude time-periodic far-field pressure.
We established the existence of a locally exponentially stable manifold of time-periodic solutions. 

Recently, we investigated the general (nonspherically symmetric) dynamics in the irrotational (less-constrained) framework \cite{LW-ipug}; 
we analyzed the linearized approximate model via spherical harmonic expansions and showed that the only possible damping mechanism in the shape / multipole modes is viscous damping.
This phenomenon is due to the spatial uniform gas pressure approximation.
Moreover, we proved the ill-posedness result arising from the viscosity-irrotationality incompatibility.
We further proved the linear well-posedness of the approximate model in the inviscid case.

Uniqueness and the characterization of equilibrium solutions of the approximate model have not been extensively studied.  
To the best of our knowledge, the only results in this direction are those by the authors: 
\begin{quote}
{\it
The equilibrium bubble must be spherical if the liquid viscosity is present \cite[Proposition 4.3\,(1)]{LW-vbas2022} or if the liquid flow is irrotational \cite[Theorem 3.1]{LW-ipug}.
}
\end{quote}
While previous studies have focused on spherically symmetric equilibrium solutions, there is a lack of understanding regarding nonspherically symmetric solutions in rotating fluids.

The main goal of this paper is to explore nonspherically symmetric equilibrium solutions of the approximate model.
Given our previous results, we neglect liquid viscosity and work with the general (rotational) frame work; see Section \ref{sec-ss-model} below.

\subsection{The steady-state system of the inviscid approximate model}\label{sec-ss-model}

Consider the steady-state system of \cite[(3.1)--(3.4)]{LW-vbas2022} (equivalently \cite[(2.5)--(2.9)]{LW-ipug}) in the inviscid case ($\mu_l=0$):

 \begin{subequations}\label{eq1.1simplified-red-ss}
\begin{empheq}[right=\empheqrbrace\text{in $\R^3\setminus \Om_*$,}]{align}
{\bf v}_{l,*}\cdot\nb{\bf v}_{l,*} =&\,- \dfrac1{\rho_l}\, \nb p_{l,*}, \label{eq1.1simplified-red-ss-a}\\
\div {\bf v}_{l,*} =&\, 0, \label{eq1.1simplified-red-ss-b}
\end{empheq}
\end{subequations}

\begin{subequations}\label{eq1.2simplified-ss-red}
\begin{empheq}[right=\empheqrbrace\text{in $\Om_*$,}]{align}
&\div(\rho_{g,*}{\bf v}_{g,*}) = 0,\qquad p_{g,*}(x) = p_{g,*},\ \text{ a constant}\label{eq1.2simplified-red-ss-a}\\
&0 = \frac{\ka}{\ga c_v} \De \log\rho_{g,*} - \frac{\ka}{\ga c_v} \frac{|\nb\rho_{g,*}|^2}{\rho_{g,*}^2} - {\bf v}_{g,*}\cdot\nb\rho_{g,*}, \label{eq1.2simplified-red-ss-c}
\end{empheq}
\end{subequations}
and
 \begin{subequations}\label{eq1.3simplified-ss-red}
\begin{empheq}[right=\empheqrbrace\text{on $\pd\Om_*$,}]{align}
{\bf v}_{l,*}\cdot\hat{\bf n} = {\bf v}_{g,*}\cdot\hat{\bf n} = 0, \label{eq1.3simplified-red-ss-a}\\
p_{g,*} - p_{l,*} = \si (\nb_S\cdot \hat {\bf n}), \label{eq1.3simplified-red-ss-b}\\
p_{g,*} = \Rg T_\infty\rho_{g,*}, \label{eq1.3simplified-red-ss-c}
\end{empheq}
\end{subequations}
with the far-field conditions 
\EQ{\label{eq-far-field-all-ss}
\lim_{|x|\to\infty} {\bf v}_{l,*}(x) = O(|x|^{-2}),\qquad
\lim_{|x|\to\infty} \nb{\bf v}_{l,*}(x) = \mathbb O,\qquad
\lim_{|x|\to\infty} p_{l,*}(x) = p_{\infty,*},
}
where the unknowns are the equilibrium bubble region $\Om_*$ and the state variables:
${\bf v}_{l,*}$, the equilibrium liquid velocity,
$p_{l,*}$, the equilibrium liquid pressure,
$\rho_{g,*}$, the equilibrium gas density, 
${\bf v}_{g,*}$, the equilibrium gas velocity, and
$p_{g,*}$, the equilibrium gas pressure.
Here, $\hat{\bf n}$ is the unit outer normal on $\Om_*$, and
$\nb_S\cdot$ denotes the surface divergence so that $\nb_S\cdot\hat{\bf n}$ is twice the mean curvature on the surface.

The system depends on physical parameters:
the density of the liquid $\rho_l>0$,
the thermal conductivity of the gas $\ka\ge0$,
the surface tension $\si$, 
far-field temperature $T_\infty>0$,
far-field pressure $p_{\infty,*}>0$,
the specific gas constant $\Rg>0$,
the heat capacity of the gas at constant volume $c_v>0$,
and the adiabatic constant $\ga>1$,
where $\Rg$, $c_v$, and $\ga$ are related by $\ga = 1 + \frac{\Rg}{c_v}$.
 
In view of the conservation of mass \cite[Proposition 7.3]{LW-vbas2022} for the evolutionary system, the steady-state system \eqref{eq1.1simplified-red-ss}--\eqref{eq-far-field-all-ss} should be parametrized by the bubble mass $M$:
\EQ{\label{eq-mass-converve}
\int_{\Om_*} \rho_{g,*}\, dx = M > 0.
}

\subsection{Overview of results}
{\ }\smallskip
Our results can be summarized as follows:

(I) \emph{Existence of an equilibrium solution with nontrivial liquid flow.}
Building upon the example of compactly supported steady incompressible Euler flows given by Gavrilov \cite{Gavrilov-GAFA2019}, we construct a solution of \eqref{eq1.1simplified-red-ss}--\eqref{eq-far-field-all-ss} with nontrivial equilibrium liquid flow ${\bf v}_{l,*}\not\equiv{\bf 0}$ and spherical bubble $\Om_* = B_{R_*}$ (Theorem \ref{prop-equilib-nt}).

(II) \emph{Characterization of spherical equilibrium bubbles in an inviscid rotational fluid.}
Previous results (\cite[Proposition 4.3\,(1)]{LW-vbas2022} and \cite[Theorem 3.1]{LW-ipug}), which characterize equilibrium bubbles as spherical, assume either $\mu_l>0$ or irrotationality of equilibrium liquid flow ${\bf v}_{l,*}$ and the far-field conditions \eqref{eq-far-field-all-ss}.
In Part (a) of Theorem \ref{thm-main}, we extend these results to axisymmetric, purely azimuthal flows, i.e., ${\bf v}_{l,*} = v_\varphi\,\hat{\boldsymbol{\varphi}}$, where $v_\varphi = v_\varphi(r,\th)$.
We show that the equilibrium bubble $\Om_*$ is spherical, provided the steady-state solution $({\bf v}_{l,*}, p_{l,*})$ belongs to the class $C^1$, even when $\mu_l=0$ and $\curl {\bf v}_{l,*}\neq{\bf 0}$ (see Remark \ref{rmk:rotation}).

(III) \emph{Construction of a horn-torus-shaped equilibrium bubble under relaxed spatial decay conditions.}
When the liquid flow has mild spatial decay (i.e.~when the spatial decay rate assumption $\lim_{|x|\to\infty} {\bf v}_{l,*}(x) = O(|x|^{-2})$ is relaxed to $\lim_{|x|\to\infty} {\bf v}_{l,*}(x) = {\bf 0}$), we construct a weak solution of \eqref{eq1.1simplified-red-ss}--\eqref{eq1.3simplified-ss-red} where the equilibrium bubble $\Om_*$ is nonspherical (horn-torus-shaped); see Part (b) of Theorem \ref{thm-main}.

\subsection{Some future directions and open problems}
{\ }\smallskip

1. \emph{Impact of spatial decay on the characterization of equilibrium bubbles.}
The equilibrium liquid flows in \cite[Proposition 3.1, Theorem 4.1]{LW-ipug} and in \cite[Proposition 4.3\,(1)]{LW-vbas2022} exhibit the spatial decay ${\bf v}_{l,*}(x) = O(|x|^{-2})$.
In this case, we have spherical equilibrium bubble $\Om_*$.
In Theorem \ref{thm-main}, on the other hand, while $\Om_*$ is not spherical, ${\bf v}_{l,*}$ does not satisfy the spatial decay ${\bf v}_{l,*}(x) = O(|x|^{-2})$.

It is natural to ask whether the spatial decay ${\bf v}_{l,*}(x) = O(|x|^{-2})$, together with surface tension $\si$, would force the equilibrium bubble $\Om_*$ spherical.

\medskip
2. \emph{Ring-torus-shaped equilibrium bubbles.}
The constructed horn-torus-shaped equilibrium bubble suggests the existence torus-shaped equilibrium bubbles.
Naturally, this leads to the question of whether ring-torus-shaped equilibrium bubbles also exist. 
To explore this direction, it might be useful to use toroidal and poloidal coordinates instead of the spherical coordinates.

\medskip
3. \emph{Using neural networks to explore new classes of equilibrium bubbles.}
In the stress balance equation \eqref{eq-SBE-appendix}, there are two unknown functions: the bubble surface $R(\th)$ and the pressure fluctuation $g(s)$.
While this paper focused on fixing $g(s) = -\si/s$ to learn $R(\th)$, it is conceivable that other nontrivial pairs $(R,g)$ yield new nonspherical equilibrium bubbles. A promising direction for future study is to simultaneously approximate both functions using neural networks--e.g., $R_{\rm NN}(\th;w_R)$ and $g_{\rm NN}(s;w_g)$--within a PINN framework.
Such an approach could open the door to discovering entirely new families of equilibrium configurations that are analytically intractable.


\bigskip\noindent{\bf Acknowledgements.}
CL and MIW are supported in part by the Simons Foundation Math + X Investigator Award \#376319 (MIW). 
CL is also supported by AMS-Simons Travel Grant.
MIW is also supported in part by National Science Foundation Grant DMS-1908657 and DMS-1937254.

\section{A spherically equilibrium bubble sitting in a steadily rotating fluid}

This section is devoted to establish the existence of nontrivial equilibrium solutions based on Gavrilov's example \cite{Gavrilov-GAFA2019}.

\begin{theorem}\label{prop-equilib-nt}
There exists a solution to the steady-state system \eqref{eq1.1simplified-red-ss}--\eqref{eq-far-field-all-ss} with ${\bf v}_{l,*}\not\equiv{\bf 0}$ and $\Om_* = B_{R_*}$ for some $R_* = R_*[M]>0$, where $M>0$ is given by \eqref{eq-mass-converve}.
\end{theorem}

\begin{proof}

The proof is built on a result  of Gavrilov \cite{Gavrilov-GAFA2019}.
Let $\mathcal{C}$ be a given circle. By \cite{Gavrilov-GAFA2019}, there exists a nontrivial smooth steady-state solution ${\bf v}_{l,*}$ to the incompressible Euler equation 
\eqref{eq1.1simplified-red-ss}, supported in a ring-like neighborhood of the circle $\mathcal{C}$.
We seek an  equilibrium spherical bubble $\Om_* = B_{R_*}$, where the bubble radius $R_* = R_*[M]$  to be determined, and so that the bubble is disjoint from the support of ${\bf v}_{l,*}$.

Outside of a neighborhood of the circle $\mathcal{C}$, ${\bf v}_{l,*} = {\bf 0}$ and $p_{l,*} = p_{\infty,*}$.
Hence  $({\bf v}_{l,*},p_{l,*})$ clearly satisfies the far-field conditions \eqref{eq-far-field-all-ss} and the steady-state kinematic boundary condition ${\bf v}_{l,*}\cdot\hat{\bf n} = 0$ on $\pd B_{R_*}$.
Moreover, the stress balance equation \eqref{eq1.3simplified-red-ss-b} reduces to $p_{g,*} = p_{\infty,*} + \frac{2\si}{R_*}$.

Finally, we simply choose a constant  density $\rho_{g,*} = \frac{p_{g,*}}{\Rg T_\infty} = \frac1{\Rg T_\infty}\bke{p_{\infty,*} + \frac{2\si}{R_*}}$ and ${\bf v}_{g,*} = {\bf 0}$ to satisfy the equilibrium equations for  the gas.
By conservation of mass \eqref{eq-mass-converve}, $M = \int_{\Om_*}\rho_{g,*}\, dx = \frac1{\Rg T_\infty}\bke{p_{\infty,*} + \frac{2\si}{R_*}} \frac{4\pi}3 R_*^3$.  The equilibrium bubble radius $R_* = R_*[M]$ is then the unique positive solution to the cubic equation $p_{\infty,*}R_*^3 + 2\si R_*^2 - \frac{3\Rg T_\infty M}{4\pi} = 0$. This completes the proof. 
\end{proof}

\section{Existence of nonspherical equilibrium bubbles}\label{appx-old-model}
 
 In this section we prove that if we relax the admit liquid velocity fields,  ${\bf v}_{l,*}$, which are slowly decaying at infinity,
  then, in the approximate model, nonspherically symmetric equilibrium bubbles are possible.
Specifically, we relax spatial decay rate assumption of ${\bf v}_{l,*}$ in \eqref{eq-far-field-all-ss}  to
\EQ{\label{eq-far-field-all}
\lim_{|x|\to\infty} {\bf v}_{l,*}(x) = {\bf 0},\qquad
\lim_{|x|\to\infty} p_{l,*}(x) = p_{\infty,*}>0.
}

\begin{theorem}\label{thm-main}
(a) Suppose that ${\bf v}_{l,*} = v_\varphi\,\hat{\boldsymbol{\varphi}}$, where $v_\varphi = v_\varphi(r,\th)$, and that the convergence of $p_{l,*}$ in \eqref{eq-far-field-all} holds.
Assume $\sigma\neq0$.
If ${\bf v}_{l,*}$ and $p_{l,*}$ are of class $C^1$, then ${\bf v}_{l,*} = {\bf 0}$, $p_{l,*} = p_{\infty,*}$ and $\Om_*$ is a sphere.

(b)
There exists a nonspherically symmetric weak solution to the steady-state system \eqref{eq1.1simplified-red-ss}--\eqref{eq1.3simplified-ss-red}.
More precisely, introduce spherical coordinates $(r,\th,\varphi)$ where $\th$ denotes the polar angle  and $\varphi$ denotes the azimuthal angle. 
Then, the following steady-state is a weak solution of  \eqref{eq1.1simplified-red-ss}--\eqref{eq1.3simplified-ss-red}:
\begin{subequations}
\label{eq-equilibrium}
\begin{align}
{\bf v}_{l,*} &= \sqrt{\frac{\si}{\rho_l r\sin\th}}\, \hat{\boldsymbol{\varphi}},\qquad\qquad\quad \ 
p_{l,*} = p_{\infty,*} - \frac{\si}{r\sin\th},\qquad\quad\
\pd\Om_*: r = C \sin\th, \label{eq-equilibrium-a}\\
\rho_{g,*} &= \frac1{\Rg T_\infty} \bke{p_{\infty,*} - \frac{4\si}C},\qquad 
{\bf v}_{g,*} = {\bf 0},\qquad\qquad\qquad\quad\quad \ \
p_{g,*} = p_{\infty,*} - \frac{4\si}{C}. \label{eq-equilibrium-b}
\end{align}
\end{subequations}
Here,  $M$ is the bubble mass given in \eqref{eq-mass-converve} and $C = C[M]>0$ denotes the unique solution of the cubic equation \eqref{eq-cubic} (below)
which is greater than $\frac{4\si}{p_{\infty,*}}$. 
 The weak solution \eqref{eq-equilibrium-a}--\eqref{eq-equilibrium-b} satisfies the relaxed far-field condition \eqref{eq-far-field-all} in every directions other than the $x_3$-direction.

\end{theorem}


\begin{remark}
The liquid flow ${\bf v}_{l,*}$ in \eqref{eq-equilibrium-a} does not satisfy the spatial decay ${\bf v}_{l,*} = O(|x|^{-2})$ in \cite[Proposition 4.3\,(1)]{LW-vbas2022} and in \cite[(2.9)]{LW-ipug}.
\end{remark}

\begin{remark}
The lower bound $\frac{2\si}{p_{\infty,*}}$ for the radius of the horn torus $\pd\Om_*$, given by $r = C\sin\theta$, is independent of the bubble mass $M$. Interestingly, even for $M = 0$, the cubic equation \eqref{eq-cubic} admits a positive solution $C = \frac{4\si}{p_{\infty,*}}$, corresponding to a horn-torus bubble with nonzero volume but zero mass--a formal ``vacuum bubble.'' Moreover, for certain negative values of $M$, the equation still yields positive solutions for $C$, suggesting the existence of horn-torus-shaped bubbles with positive volume but negative mass.

This behavior stands in contrast to the spherical case, where the corresponding cubic equation
\[
p_\infty R_*^3 + 2\si R_*^2 - \frac{3R_g T_\infty M}{4\pi} = 0
\]
admits no positive solution $R_*$ when $M \le 0$, implying that spherical equilibrium bubbles cannot exist with zero or negative mass.

This raises the question of whether such horn-torus-shaped bubbles with $ M \le 0$ represent physically realizable states--such as vacuum bubbles--or if they are artifacts of the asymptotic model in extreme regimes. Further investigation is needed to determine whether the mass-independent lower bound and the existence of bubbles with negative mass reflect real physical phenomena or limitations of the current model.
\end{remark}

\begin{remark}\label{rmk:rotation}
In the proof of Theorem \ref{thm-main}, ${\bf v}_{l,*}$ is not necessarily irrotational. 
In fact, for ${\bf v}_{l,*} = v_\varphi\,\hat{\boldsymbol{\varphi}}$, where $v_\varphi = v_\varphi(r,\th)$, we have
\EQN{
\curl{\bf v}_{l,*} &= 
\frac1{r\sin\th} \pd_\th(v_\varphi\sin\th)\,\hat{\bf r} - \frac1r \pd_r(rv_\varphi)\,\hat{\boldsymbol{\th}}\\
&= \bke{ \frac{\pd_\th v_\varphi}{r} + \frac{\cot\th}{r} v_\varphi} \hat{\bf r} - \bke{\pd_rv_\varphi + \frac{v_\varphi}{r}} \hat{\boldsymbol{\th}}
\neq {\bf 0},\ \ \text{in general.}
}
For example, for $v_\varphi(r,\th) = \frac1{r}$, 
\[
\curl{\bf v}_{l,*} = \frac{\cot\th}{r^2} \,\hat{\bf r} \neq {\bf 0}.
\]
Moreover, the weak solution \eqref{eq-equilibrium} is not irrotational: since $v_\varphi = \sqrt{\frac{\si}{\rho_l r\sin\th}}$,
\[
\curl{\bf v}_{l,*} = \frac12\sqrt{\frac{\si}{\rho_l}}\frac{\cot\th}{r^{\frac32}\sqrt{\sin\th}}\, \hat{\bf r} + \frac12\sqrt{\frac{\si}{\rho_l}}\frac1{r^{\frac32}\sqrt{\sin\th}}\, \hat{\boldsymbol{\th}} 
\neq {\bf 0}.
\]
\end{remark}

\begin{proof}[Proof of Theorem \ref{thm-main}]
Using spherical coordinates ${\bf v}_{l,*} = v_r\,\hat{\bf r} + v_\th\,\hat{\boldsymbol{\th}} + v_\varphi\,\hat{\boldsymbol{\varphi}}$, $\th\in(0,\pi)$, $\varphi\in(0,2\pi)$,
the incompressible Euler equation for the liquid reads
 \begin{subequations}\label{eq-euler-spherical}
\begin{empheq}{align}
v_r\pd_rv_r + \frac{v_\th}r\pd_\th v_r + \frac{v_\varphi}{r\sin\th}\pd_\varphi v_r - \frac{v_\th^2+v_\varphi^2}r &= -\frac1{\rho_l}\pd_r p_{l,*},\\
v_r\pd_rv_\th + \frac{v_\th}r\pd_\th v_\th + \frac{v_\varphi}{r\sin\th}\pd_\varphi v_\th + \frac{v_\th v_\varphi}r - \frac{v_\varphi^2\cot\th}r &= -\frac1{\rho_l r} \pd_\th p_{l,*},\\
v_r\pd_r v_\varphi + \frac{v_\th}r \pd_\th v_\varphi + \frac{v_\varphi}{r\sin\th}\pd_\varphi v_\varphi + \frac{v_rv_\varphi}r + \frac{v_\th v_\varphi\cot\th}r &= -\frac1{\rho_l r\sin\th}\pd_\varphi p_{l,*},\\
\frac1{r^2}\pd_r(r^2v_r) + \frac1{r\sin\th}\pd_\th(\sin\th v_\th) + \frac1{r\sin\th}\pd_\varphi v_\varphi &= 0.
\end{empheq}
\end{subequations}
Note that for $C^1$ solutions, the regularity condition $v_\varphi\sim\sin\theta$ must hold as $\theta\to0$ or $\theta\to\pi$.

Consider the ansatz $v_r=v_\th=0$, and $v_\varphi = v_\varphi(r,\th)$, $p_{l,*}=p_{l,*}(r,\th)$. 
Then \eqref{eq-euler-spherical} is reduced to 
\begin{subequations}
\begin{empheq}{align}
- \frac{v_\varphi^2}r &= -\frac1{\rho_l}\pd_r p_{l,*},\label{eq-v_varphi2-a}\\
- \frac{v_\varphi^2\cot\th}r &= -\frac1{\rho_l r} \pd_\th p_{l,*},\label{eq-v_varphi2-b}
\end{empheq}
\end{subequations}
which implies that 
\[
r\pd_rp_{l,*}\cot\th = \pd_\th p_{l,*}.
\]
By the method of characteristics, we obtain
\EQ{\label{eq-pl}
p_{l,*}(r,\th) = p_{\infty,*} + g(r\sin\th),\quad \textrm{with\ \  $g\in C^1((0,\infty))$\ \  and \ \ $\lim_{s\to\infty}g(s) = 0$, }
}
to be determined.
Here, we have used the convergence of $p_{l,*}$ in \eqref{eq-far-field-all}.
Note that \eqref{eq-pl} holds for $\theta\in(0,\pi)$.
From \eqref{eq-pl}, it follows that
\EQ{\label{eq-vl}
{\bf v}_{l,*} = \sqrt{\frac{r g'(r\sin\th)\sin\th}{\rho_l}}\, \hat{\boldsymbol{\varphi}},\qquad 
\th\in(0,\pi).
}
If the solution is of class $C^1$, \eqref{eq-pl} and \eqref{eq-vl} extend to $\theta\in[0,\pi]$, in which case $g\in C^1([0,\infty))$.

Note that, using \eqref{eq-pl} in \eqref{eq-v_varphi2-a}, or in \eqref{eq-v_varphi2-b} alternatively, yields $g'(r\sin\th) = \frac{\rho_l v_\varphi^2}{r\sin\th} \ge0$ for $\th\in(0,\pi)$, which implies $g'(s)\ge0$ for $s>0$.
Since $\lim_{s\to\infty}g(s)=0$ and $g'(s)\ge0$, we have $g(s)\le0$ for all $s>0$.

Suppose the solution is $C^1$, then using the regularity condition $v_\varphi\sim\sin\th$, as $\theta\to0$ or $\theta\to\pi$, in \eqref{eq-v_varphi2-b} yields $\pd_rp_{l,*}(r,0) = 0$, which implies that $p_{l,*}(r,0) = p_{\infty,*}$ for all $r$.
Hence, by evaluating \eqref{eq-pl} at $\th=0$, we obtain $g(0)=0$, which, combining with the fact that $g'(s)\ge0$ and $\lim_{s\to\infty}g(s)=0$, gives $g\equiv0$. 
Therefore, we have $p_{l,*}(r,\th) \equiv p_{\infty,*}$ by using $g\equiv0$ in \eqref{eq-pl}.
Thus by Laplace--Young condition \eqref{eq1.3simplified-red-ss-b}, 
\[
(\nb_S\cdot\hat{\bf n})(\th) = \sigma^{-1} \left(p_{g,*}-p_{l,*}(r,\th)\Big|_{r=R(\th)}\right) = \sigma^{-1} (p_{g,*}-p_{\infty,*}) = 
\text{constant},
\]
implying that $\Om_*$ is spherical by Alexandrov's theorem \cite{Alexandrov-AMPA1962}.
This proves Part (a) of the theorem.

Let the equilibrium bubble surface $\pd\Om_*$ be parametrized by $F(r,\th,\varphi) = R(\th,\varphi) - r = 0$.
If we make the ansatz $R(\th,\varphi) = R(\th)$, then the unit normal on $\pd\Om_*$ is given by
\[
\hat{\bf n} = \frac{\nb F}{|\nb F|} = \frac{ \pd_rF\, \hat{\bf r} + \frac1r\pd_\th F\, \hat{\boldsymbol{\th}} + \frac1{r\sin\th}\pd_\varphi F\, \hat{\boldsymbol{\varphi}}}{|\nb F|}
= \frac{- \hat{\bf r} + \frac1rR'(\th)\, \hat{\boldsymbol{\th}} }{|\nb F|},
\]
so that ${\bf v}_{l,*}\cdot\hat{\bf n} = 0$ on $\pd\Om_*$.
Here we've extended the domain of $\hat{\bf n}$ to a neighborhood of $\pd\Om_*$ so that $\hat{\bf n} = \hat{\bf n}(r,\th,\varphi)$.

Write 
\[
\hat{\bf n} = n_r\,\hat{\bf r} + n_\th\,\hat{\boldsymbol{\th}} + n_\varphi\,\hat{\boldsymbol{\varphi}} 
= \frac{-1}{\sqrt{1+\frac{(R'(\th))^2}{r^2}}}\, \hat{\bf r} + \frac{\frac{R'(\th)}r}{\sqrt{1+\frac{(R'(\th))^2}{r^2}}}\, \hat{\boldsymbol{\th}} + 0\,\hat{\boldsymbol{\varphi}}.
\]

Next  we compute the mean curvature term $\nb_S\cdot\hat{\bf n}$.   For the choice
\[  R(\th) = C\sin\th,\ 0<\theta <\pi, \]
where  $C>0$  is a constant, we obtain the relation 
\begin{equation}
\nb_S\cdot\hat{\bf n} = \frac1C\bke{\frac1{\sin^2\th} - 4}.
\label{eq:div-n}
\end{equation}
Two versions of the computation are presented in Appendix \ref{appendix-curvature}.

Note that the surface $\pd\Om_*: r=C\sin\th$ is a horn torus with both major and minor radii being equal to $C/2$ (see e.g. \cite[pp. 305--306]{horn-torus-book}). 

We next impose Laplace--Young condition \eqref{eq1.3simplified-red-ss-b}.
 For the left hand side of \eqref{eq1.3simplified-red-ss-b} we have, using \eqref{eq-pl}:
\[ p_{g,*}-p_{l,*}  = p_{g,*}-p_{\infty,*}  - g\big(r\sin\theta\big)\Big|_{r=C\sin\theta}\  =\  p_{g,*}-p_{\infty,*}  - g\big(C\sin^2\theta\big)\]
and for the right hand side of \eqref{eq1.3simplified-red-ss-b} we have, using the relation \eqref{eq:div-n}:
\[ \sigma \nb_S\cdot\hat{\bf n} = \frac{\sigma}{C}\bke{\frac1{\sin^2\th} - 4}. \]
Hence, \eqref{eq1.3simplified-red-ss-b} holds if 
\[ p_{g,*}-p_{\infty,*}  - g\big(C\sin^2\theta\big) = \frac{\sigma}{C}\bke{\frac1{\sin^2\th} - 4}\ .\]
For any given $p_{\infty,*}$, we choose $p_{g,*} = p_{\infty,*} - 4\sigma/C$. This imposes the following condition on $g$:
\[ -g(C\sin^2\theta) = \frac{\sigma}{C\sin^2\theta}, \]
which we satisfy by setting
\[ g(s) \equiv - \frac{\sigma}{s}.\]

Further, we may take ${\bf v}_{g,*} = {\bf 0}$ and, by the ideal gas law: 
 \[ \rho_{g,*}\equiv\frac{p_{g,*}}{\Rg T_\infty} = \frac1{\Rg T_\infty} \bke{p_{\infty,*} - \frac{4\si}C} \]

It remains to determine the positive constant, $C$.
Note that the volume of $|\Om_*|$, when the boundary is parametrized by $r=R(\th)$, is given by
\[
|\Om_*| = \int_0^{2\pi}\int_0^\pi \int_0^{R(\th)} r^2\sin\th\, drd\th d\varphi\\
= \frac{2\pi}3 \int_0^\pi (R(\th))^3\sin\th\, d\th.
\]
So, if $R(\th) = C\sin\th$, we have
\[
|\Om_*| = \frac{2\pi}3 C^3 \int_0^\pi \sin^4\th\, d\th = \frac{\pi^2}4 C^3.
\]
 By conservation of mass formula \eqref{eq-mass-converve},
\EQ{\label{eq-mass-formula}
M = \rho_{g,*}|\Om_*| = \rho_{g,*} \frac{\pi^2}4 C^3  =\frac{p_{\infty,*} - \frac{4\si}C}{\Rg T_\infty}\, \frac{\pi^2}{4}  C^3.
}
Note that we require $C>\frac{4\si}{p_{\infty,*}}$ to ensure $M>0$.
Moreover, \eqref{eq-mass-formula} is equivalent to the cubic equation
\EQ{\label{eq-cubic}
p_{\infty,*} C^3 - 4\si C^2 - \frac{4\Rg T_\infty M}{\pi^2} = 0.
}
Let $f(C) = p_{\infty,*} C^3 - 4\si C^2 - \frac{4\Rg T_\infty M}{\pi^2}$. 
Since the leading coefficient $p_{\infty,*} > 0$, we have $f(C)>0$ for sufficiently large $C>0$.
Evaluating $f$ at $\frac{4\si}{p_{\infty,*}}$, we find $f(\frac{4\si}{p_{\infty,*}}) = - \frac{4\Rg T_\infty M}{\pi^2} < 0$.
Therefore, by the intermediate value theorem (Bolzano's theorem), there exists a root of $f(C) = 0$ in the interval $(\frac{4\si}{p_{\infty,*}}, +\infty)$.
This solution is unique. Indeed, the derivative of $f$ is given by
\[
f'(C) = 3p_{\infty,*} C^2 - 8\si C,
\]
which vanishes only at $C = 0$ and $C = \frac{8\si}{3p_{\infty,*}}$.
Since $f(C)$ is a cubic polynomial with positive leading coefficient, it has only one local maximum at the smaller critical point $C=0$, where $f_{\rm local max} = f(0) = - \frac{4\Rg T_\infty M}{\pi^2}<0$.
Hence, $f(C)=0$ has exactly one root in $(\frac{4\si}{p_{\infty,*}}, +\infty)$.
We choose $C$ to be the unique solution of \eqref{eq-cubic} that satisfies $C>\frac{4\si}{p_{\infty,*}}$.

We now prove that \eqref{eq-equilibrium-a}--\eqref{eq-equilibrium-b} is a weak solution of \eqref{eq1.1simplified-red-ss}--\eqref{eq1.3simplified-ss-red} in the sense that 
${\bf v}_{l,*}\in L^2_{\rm loc}(\R^3\setminus\Om_*)$, meets the boundary condition \eqref{eq1.3simplified-red-ss-a} pointwise, and satisfies the weak form
\[
\int_{\R^3\setminus\Om_*} {\bf v}_{l,*}^\top\cdot\nb\boldsymbol{\ze}\cdot{\bf v}_{l,*}\, dx = 0,\qquad \forall \boldsymbol{\ze}\in {\bf C}^\infty_c(\R^3\setminus\overline{\Om_*})\text{ with }\div\boldsymbol{\ze} = 0,
\]
\[
\int_{\R^3\setminus\Om_*} {\bf v}_{l,*}\cdot\nb\phi\, dx = 0,\qquad \forall \phi\in C^\infty_c({\R^3\setminus\overline{\Om_*}}).
\]
Indeed, the functions in \eqref{eq-equilibrium} satisfy the system in the classical sense everywhere except on the $x_3$-axis since ${\bf v}_{l,*}$ is singular along the $x_3$-axis.
Nontheless, ${\bf v}_{l,*}$ is $L^2_{\rm loc}$ since the Jacobian is $r^2 \sin\th$, and so it satisfies the weak formulation.
More precisely, if $\boldsymbol{\ze} = \ze_r\,\hat{\bf r} + \ze_\th\,\hat{\boldsymbol{\th}} + \ze_\varphi\,\hat{\boldsymbol{\varphi}}$, then 
\EQN{
{\bf v}_{l,*}^\top\cdot\nb\boldsymbol{\ze}\cdot{\bf v}_{l,*} 
&= \frac{\si}{\rho_l r\sin\th} \bke{\frac1{r\sin\th}\pd_\varphi\ze_\varphi + \cot\th\frac{\ze_\th}r + \frac{\ze_r}r}\\
&= \frac{\si}{\rho_l r\sin\th} \bke{-\pd_r\ze_r - \frac1r\ze_r - \frac1r\pd_\th\ze_\th},
}
where we've used $0 = \div\boldsymbol{\ze} = \pd_r\ze_r + \frac2r\ze_r + \cot\th\frac{\ze_\th}r + \frac1r \pd_\th\ze_\th + \frac1{r\sin\th}\pd_\varphi\ze_\varphi$.
Therefore, we have
\EQN{
\int_{\R^3\setminus\Om_*} {\bf v}_{l,*}^\top\cdot\nb\boldsymbol{\ze}\cdot{\bf v}_{l,*}\, dx
&= \iiint_{\text{supp}(\boldsymbol{\ze})} \frac{\si}{\rho_l r\sin\th} \bke{-\pd_r\ze_r - \frac1r\ze_r - \frac1r\pd_\th\ze_\th} r^2\sin\th\, drd\th d\varphi\\
&= \frac{\si}{\rho_l} \iiint_{\text{supp}(\boldsymbol{\ze})} \bkt{-\pd_r\bke{r\ze_r} - \pd_\th\ze_\th}  drd\th d\varphi = 0,
}
by using the integration by parts formula and compactness of the support of $\boldsymbol{\ze}$.

Moreover, since $\nb\phi = \pd_r\phi\,\hat{\bf r} + \frac1r\pd_\th\phi\,\hat{\boldsymbol{\th}} + \frac1{r\sin\th} \pd_\varphi\phi\,\hat{\boldsymbol{\varphi}}$,
\EQN{
\int_{\R^3\setminus\Om_*}  {\bf v}_{l,*}\cdot\nb\phi\, dx 
&= \iiint_{\text{supp}(\phi)} \sqrt{\frac{\si}{\rho_l r\sin\th}} \bke{\frac1{r\sin\th}\pd_\varphi\phi} r^2\sin\th\, dr d\th d\varphi\\
&= \sqrt{\frac{\si}{\rho_l}} \iiint_{\text{supp}(\phi)} \frac{\pd_\varphi\phi}{\sqrt{\sin\th}} \sqrt{r}\, dr d\th d\varphi\\
&= \sqrt{\frac{\si}{\rho_l}} \iint \frac{\sqrt{r}}{\sqrt{\sin\th}} \bke{\int_0^{2\pi} \pd_\varphi\phi\, d\varphi} dr d\th = 0,
}
by using the integration by parts formula and compactness of the support of $\phi$.

This completes the proof of Theorem \ref{thm-main}.
\end{proof}

\begin{remark}
One could use a more general ansatz $R(\th) = C_1\cos\th + C_2\sin\th$ to compute the curvature $\nb_S\cdot\hat{\bf n}$ in \eqref{eq-curvature} and match the Laplace--Young condition.
However, the constraint $R(\th) = r \ge0$ for $\th\in[0,\pi]$ implies that $C_1=0$.
\end{remark}

\appendix

\section{Computation of the mean curvature term}\label{appendix-curvature}

\subsection{$\nabla_S\cdot \hat{\bf n}$ via extension of $\hat{\bf n}$ off the surface}

\EQ{\label{eq-curvature}
\nb_S\cdot\hat{\bf n} &=  \div\hat{\bf n} \Big|_{r=R}  = \bkt{ \frac1{r^2}\pd_r(r^2n_r) + \frac1{r\sin\th}\pd_\th(\sin\th\, n_\th) + \frac1{r\sin\th}\pd_\varphi n_\varphi } \Big|_{r=R}\\
& = - \frac{2R^2+3(R')^2}{(R^2+(R')^2)^{\frac32}} + \frac{\cos\th R'R^2 + \cos\th(R')^3 + \sin\th R^2R'' + 2\sin\th(R')^2R''}{(R^2+(R')^2)^{\frac32}R\sin\th}\\
& = \frac{- 2\sin\th R^3 - 3\sin\th R(R')^2 + \cos\th R'R^2 + \cos\th(R')^3 + \sin\th R^2R'' + 2\sin\th(R')^2R''}{(R^2+(R')^2)^{\frac32}R\sin\th}.
}

For $R(\th)$, we use the ansatz $R(\th) = C\sin\th$, where $C>0$ to be determined, and get
\[
\nb_S\cdot\hat{\bf n} = \frac1C\bke{\frac1{\sin^2\th} - 4}.
\]

\subsection{$\nabla_S\cdot \hat{\bf n}$ via the second fundamental form}

In this appendix, we provide a different computation of the surface divergence in \eqref{eq-curvature} without extending the domain of the normal $\hat{\bf n}$.
We instead use the formula of mean curvature given in terms of first and second fundamental forms.

Let the equilibrium bubble surface $\pd\Om_*$ be parametrized by $r = R(\th,\varphi)$ in spherical coordinates.
Let ${\bf x}:\Om_*\to\R^3$ be a regular patch.
Then 
\[
{\bf x} = 
\begin{bmatrix}
R(\th,\varphi)\sin\th\cos\varphi\\
R(\th,\varphi)\sin\th\sin\varphi\\
R(\th,\varphi)\cos\th
\end{bmatrix}.
\]
We compute the mean curvature $H$ by using the formula:
\[
H = \frac{eG - 2fF + gE}{2(EG - F^2)},
\]
where $E,F,G$ are coefficients of the first fundamental form, and $e,f,g$ are coefficients of the second fundamental form.

Assume $R=R(\th)$.
Then the partial derivatives of ${\bf x}$ are
\[
{\bf x}_\th = 
\begin{bmatrix}
(R'\sin\th + R\cos\th)\cos\varphi\\
(R'\sin\th + R\cos\th)\sin\varphi\\
R'\cos\th - R\sin\th
\end{bmatrix},\qquad
{\bf x}_\varphi = 
\begin{bmatrix}
-R\sin\th\sin\varphi\\
R\sin\th\cos\varphi\\
0
\end{bmatrix},
\]
where $R' = \frac{d}{d\th} R(\th)$.
Thus, the coefficients of the first fundamental form are
\[
E = \norm{{\bf x}_\th}^2 = (R')^2 + R^2,\qquad
F = {\bf x}_\th\cdot{\bf x}_\varphi = 0,\qquad
G = \norm{{\bf x}_\varphi}^2 = R^2\sin^2\th,
\]
and $EG - F^2 = \bkt{(R')^2+R^2}R^2\sin^2\th$.

Moreover, the second-order partial derivatives of ${\bf x}$ are
\[
{\bf x}_{\th\th} = 
\begin{bmatrix}
(R''\sin\th + 2R'\cos\th - R\sin\th)\cos\varphi\\
(R''\sin\th + 2R'\cos\th - R\sin\th)\sin\varphi\\
R''\cos\th - 2R'\sin\th - R\cos\th
\end{bmatrix},
\]
\[
{\bf x}_{\th\varphi} = 
\begin{bmatrix}
-(R'\sin\th+R\cos\th)\sin\varphi\\
(R'\sin\th+R\cos\th)\cos\varphi\\
0
\end{bmatrix},\qquad
{\bf x}_{\varphi\varphi} = 
\begin{bmatrix}
-R\sin\th\cos\varphi\\
-R\sin\th\sin\varphi\\
0
\end{bmatrix}.
\]
Recall that the coefficients of the second fundamental form are given by
\[
e = \frac{\det\bke{ {\bf x}_{\th\th}\, {\bf x}_\th\, {\bf x}_\varphi }}{\sqrt{EG-F^2}},\qquad
f = \frac{\det\bke{ {\bf x}_{\th\varphi}\, {\bf x}_\th\, {\bf x}_\varphi }}{\sqrt{EG-F^2}},\qquad
g = \frac{\det\bke{ {\bf x}_{\varphi\varphi}\, {\bf x}_\th\, {\bf x}_\varphi }}{\sqrt{EG-F^2}}.
\]
We compute
\EQN{
\det\bke{ {\bf x}_{\th\th}\, {\bf x}_\th\, {\bf x}_\varphi } 
&= 
\begin{vmatrix}
(R''\sin\th + 2R'\cos\th - R\sin\th)\cos\varphi&(R'\sin\th + R\cos\th)\cos\varphi&-R\sin\th\sin\varphi\\
(R''\sin\th + 2R'\cos\th - R\sin\th)\sin\varphi&(R'\sin\th + R\cos\th)\sin\varphi&R\sin\th\cos\varphi\\
R''\cos\th - 2R'\sin\th - R\cos\th&R'\cos\th - R\sin\th&0
\end{vmatrix}\\
&= 
\begin{vmatrix}
(R''\sin\th + 2R'\cos\th - R\sin\th)\cos\varphi&(R'\sin\th + R\cos\th)\cos\varphi&-R\sin\th\sin\varphi\\
(R''\sin\th + 2R'\cos\th - R\sin\th)\frac1{\sin\varphi}&(R'\sin\th + R\cos\th)\frac1{\sin\varphi}&0\\
R''\cos\th - 2R'\sin\th - R\cos\th&R'\cos\th - R\sin\th&0
\end{vmatrix}\\
&= -R\sin\th\sin\varphi 
\begin{vmatrix}
(R''\sin\th + 2R'\cos\th - R\sin\th)\frac1{\sin\varphi}&(R'\sin\th + R\cos\th)\frac1{\sin\varphi}\\
R''\cos\th - 2R'\sin\th - R\cos\th&R'\cos\th - R\sin\th
\end{vmatrix}\\
&= R\sin\th\bkt{R''R -2(R')^2 - R^2},
}

\EQN{
\det\bke{ {\bf x}_{\th\varphi}\, {\bf x}_\th\, {\bf x}_\varphi } 
&= 
\begin{vmatrix}
-(R'\sin\th+R\cos\th)\sin\varphi&(R'\sin\th + R\cos\th)\cos\varphi&-R\sin\th\sin\varphi\\
(R'\sin\th+R\cos\th)\cos\varphi&(R'\sin\th + R\cos\th)\sin\varphi&R\sin\th\cos\varphi\\
0&R'\cos\th - R\sin\th&0
\end{vmatrix}\\
&= -\bke{R'\cos\th - R\sin\th}
\begin{vmatrix}
-(R'\sin\th+R\cos\th)\sin\varphi&-R\sin\th\sin\varphi\\
(R'\sin\th+R\cos\th)\cos\varphi&R\sin\th\cos\varphi\\
\end{vmatrix}\\
&= 0,
}

\EQN{
\det\bke{ {\bf x}_{\varphi\varphi}\, {\bf x}_\th\, {\bf x}_\varphi } 
&= 
\begin{vmatrix}
-R\sin\th\cos\varphi&(R'\sin\th + R\cos\th)\cos\varphi&-R\sin\th\sin\varphi\\
-R\sin\th\sin\varphi&(R'\sin\th + R\cos\th)\sin\varphi&R\sin\th\cos\varphi\\
0&R'\cos\th - R\sin\th&0
\end{vmatrix}\\
&= -\bke{R'\cos\th - R\sin\th}
\begin{vmatrix}
-R\sin\th\cos\varphi&-R\sin\th\sin\varphi\\
-R\sin\th\sin\varphi&R\sin\th\cos\varphi\\
\end{vmatrix}\\
&= \bke{R'\cos\th - R\sin\th}R^2\sin^2\th.
}
Hence,
\[
e = \frac{R''R - 2(R')^2 - R^2}{\sqrt{(R')^2 + R^2}},\qquad
f = 0,\qquad
g = \frac{R\sin\th\bkt{R'\cos\th - R\sin\th}}{\sqrt{(R')^2 + R^2}}.
\]
Therefore, the mean curvature $H$ is (see eg. \cite[Eq.~(5) in p.156]{DoCarmo-book})
\EQN{
H = \frac{eG - 2fF + gE}{2(EG - F^2)}
&= \frac{\bke{R''R - 2(R')^2 - R^2}R\sin\th + \bke{R'\cos\th - R\sin\th} \bke{(R')^2 + R^2}}{2\bke{(R')^2 + R^2}^{\frac32}R\sin\th}\\
&= \frac{-2\sin\th R^3 -3\sin\th R(R')^2 + \cos\th R' R^2 + \cos\th(R')^3 + \sin\th R^2 R''}{2\bke{(R')^2 + R^2}^{\frac32}R\sin\th}.
}
This implies that
\[
\nb_S\cdot\hat{\bf n} = 2H 
= \frac{- 2\sin\th R^3 - 3\sin\th R(R')^2 + \cos\th R' R^2 + \cos\th(R')^3 + \sin\th R^2 R''}{\bke{(R')^2 + R^2}^{\frac32}R\sin\th},
\]
which is consistent with the formula \eqref{eq-curvature}.

\section{Boundary conditions for $R$}
Recall the stress balance equation (Laplace--Young condition):
\EQ{\label{eq-SBE-appendix}
p_{g,*}&-p_{\infty,*}  - g(R\sin\th) \\
&= \si\,\frac{- 2\sin\th R^3 - 3\sin\th R(R')^2 + \cos\th R'R^2 + \cos\th(R')^3 + \sin\th R^2R'' + 2\sin\th(R')^2R''}{(R^2+(R')^2)^{\frac32}R\sin\th}.
}
Assume $R, R', R''$ are bounded. 
Multiplying both sides of the equation by $-R\sin\th$ and taking the limit as $\th\to0$, we obtain 
\EQ{\label{eq-BC0}
\lim_{s\to0} g(s)s = -\frac{\si R'(0)}{\sqrt{(R(0))^2+(R'(0))^2}}.
}
Similarly, taking the limit $\th\to\pi$ yields
\EQ{\label{eq-BCpi}
\lim_{s\to0} g(s)s = \frac{\si R'(\pi)}{\sqrt{(R(\pi))^2+(R'(\pi))^2}}.
}

\section{Numerical simulation}\label{sec:numerics}

In our setting, the stress balance equation \eqref{eq-SBE-appendix} is a second-order, highly nonlinear ODE. In this appendix, we apply the Physics-Informed Neural Network (PINN) method to solve this ODE.

In general, one could attempt to learn both $R(\th)$ and $g(s)$ in \eqref{eq-SBE-appendix}. However, the target function $g(s)=-\si/s$ is singular and therefore difficult to approximate directly using neural networks. 
To simplify the problem, we fix $g(s)=-\si/s$ and focus solely on learning $R(\th)$.
Under this choice, the stress balance equation becomes
\[
(L_{\rm SB}(R))(\th) = 0,
\]
where
\EQN{
(L_{\rm SB}(R))(\th) &:= p_{g,*}-p_{\infty,*} + \frac{\si}{R\sin\th} \\
&\quad - \si\,\frac{- 2\sin\th R^3 - 3\sin\th R(R')^2 + \cos\th R'R^2 + \cos\th(R')^3 + \sin\th R^2R'' + 2\sin\th(R')^2R''}{(R^2+(R')^2)^{\frac32}R\sin\th}.
}
The boundary conditions \eqref{eq-BC0} and \eqref{eq-BCpi} reduce to
\[
\si = \frac{\si R'(0)}{\sqrt{(R(0))^2+(R'(0))^2}} = -\frac{\si R'(\pi)}{\sqrt{(R(\pi))^2+(R'(\pi))^2}},
\]
which are equivalent to
\[
R'(0) = \sqrt{(R(0))^2+(R'(0))^2},\qquad 
R'(\pi) = -\sqrt{(R(\pi))^2+(R'(\pi))^2}.
\]

Our goal is to find a function $R(\th)$ satisfying $(L_{\rm SB}(R))(\th) = 0$, subject to the prescribed bubble volume
\[
V:= \frac{M}{\rho_g} = |\Om_*| = \frac{2\pi}3 \int_0^\pi (R(\th))^3 \sin\th\, d\th =: V(R).
\]
Due to the translational invariance of the system, we impose the symmetry $R(\th) = R(\pi-\th)$ to fix the centroid at the origin.
Accordingly, we define a symmetric neural network approximation:
\[
R_{\rm NN}(\th;w) = {\rm Softplus}(L_4\circ\tanh\circ L_3\circ\tanh L_2\circ\tanh\circ L_1(\th_{\rm sym})),
\]
where 
$\th_{\rm sym} = \frac{\pi}2 - \abs{\th-\frac{\pi}2}$, and ${\rm Softplus}(x) = \log(1+\exp(x))$ is used to ensure $R_{\rm NN}(\th;w)>0$.
This symmetry allows us to solve the equation only on the interval $(0,\frac{\pi}2)$.
To ensure smoothness at $\th=\pi/2$, we enforce $R_{\rm NN}'(\frac{\pi}2;w) = 0$.
Here, $w$ denotes the collection of all weights and biases in the neural network architecture, which are the parameters to be optimized during training.

We discretize the domain by setting $\th_i = (i-1)\De\th$, where $\De\th = \frac{\pi}{2(N-1)}$ and $i=1,\ldots,N$.
The training objective is to minimize the total loss: 
\[
w^* = \underset{w}{\text{argmin}} \bke{ \la_{\rm SB}\mathcal{L}_{\rm SB}(w) + \la_V\mathcal{L}_V(w) + \la_B\mathcal{L}_B(w) + \la_S\mathcal{L}_S(w) },
\]
with the following loss terms:
\[
\mathcal{L}_{\rm SB}(w) = \frac1{N} \sum_{i=1}^N \bke{L_{\rm SB}(R_{\rm NN}(\th_i;w)) }^2,
\]
\[ 
\mathcal{L}_V(w) = \bke{ \frac{V(R_{\rm NN}(\cdot\,;w))-V}{V} }^2,\quad \textrm{in which}\ \ 
V(R_{\rm NN}(\cdot\,;w)) = \frac{2\pi}3 \De\th \sum_{i=1}^N \bke{R_{\rm NN}(\th_i;w)}^3\sin\th_i,
\] 
\[
\mathcal{L}_B(w) = \bke{ R_{\rm NN}'(0;w) - \sqrt{(R_{\rm NN}(0;w))^2 + (R_{\rm NN}(0;w))^2} }^2,
\]
and 
\[
\mathcal{L}_S(w) = \bke{ R_{\rm NN}'(\frac{\pi}2;w)}^2.
\]

To evaluate the quality of the learned solution, we compute the relative root mean square error (rRMSE) with respect to the spherical profile $R(\th) = C\sin\th$, where $C=\bke{\frac{4V}{\pi^2}}^{1/3}$:
\[
{\rm rRMSE} = \frac{\sqrt{\sum_{i=1}^N \abs{R_{\rm NN}(\th_i) - C\sin\th_i}^2}}{\sqrt{\sum_{i=1}^N \abs{C\sin\th_i}^2}}.
\]

\begin{figure}[H]
\centering
\subfloat[3D surface of the learned equilibrium bubble.]{
  \includegraphics[width=0.31\textwidth]{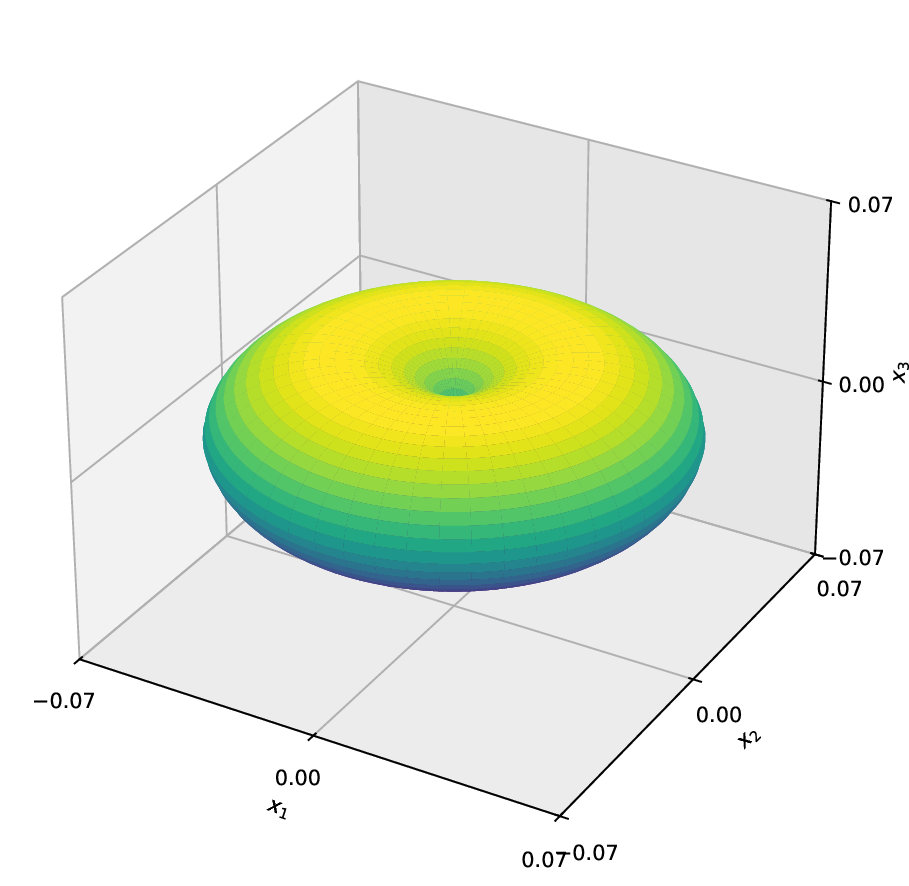}
  \label{fig:3D}
}
\hfill
\subfloat[The learned $R_{\rm NN}(\th)$ (blue curve) versus the ground truth $R(\th) = C\sin\th$ (red curve).]{
\includegraphics[width=0.31\textwidth]{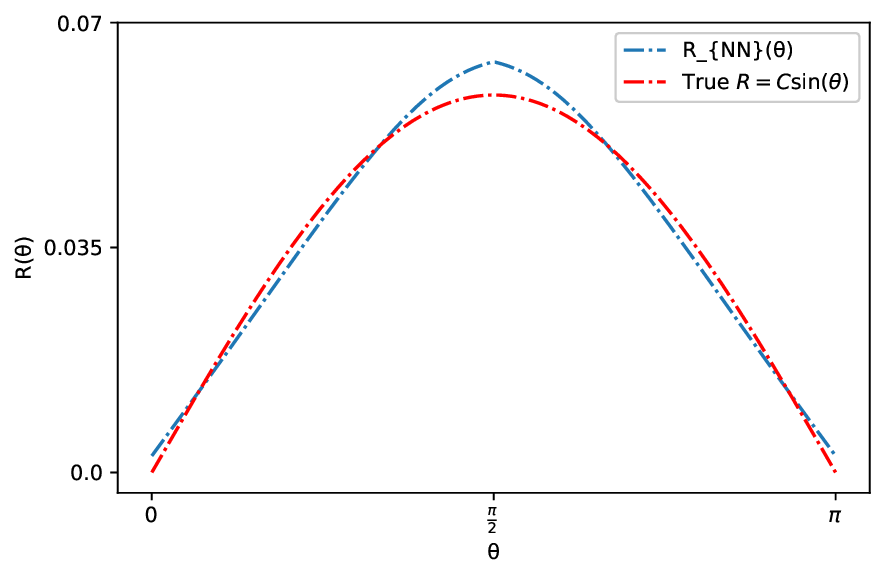}
\label{fig:R}
}
\hfill
\subfloat[Polar plot of the learned $R_{\rm NN}(\th)$ (blue curve) versus the ground truth $R(\th) = C\sin\th$ (red curve).]{
\includegraphics[width=0.31\textwidth]{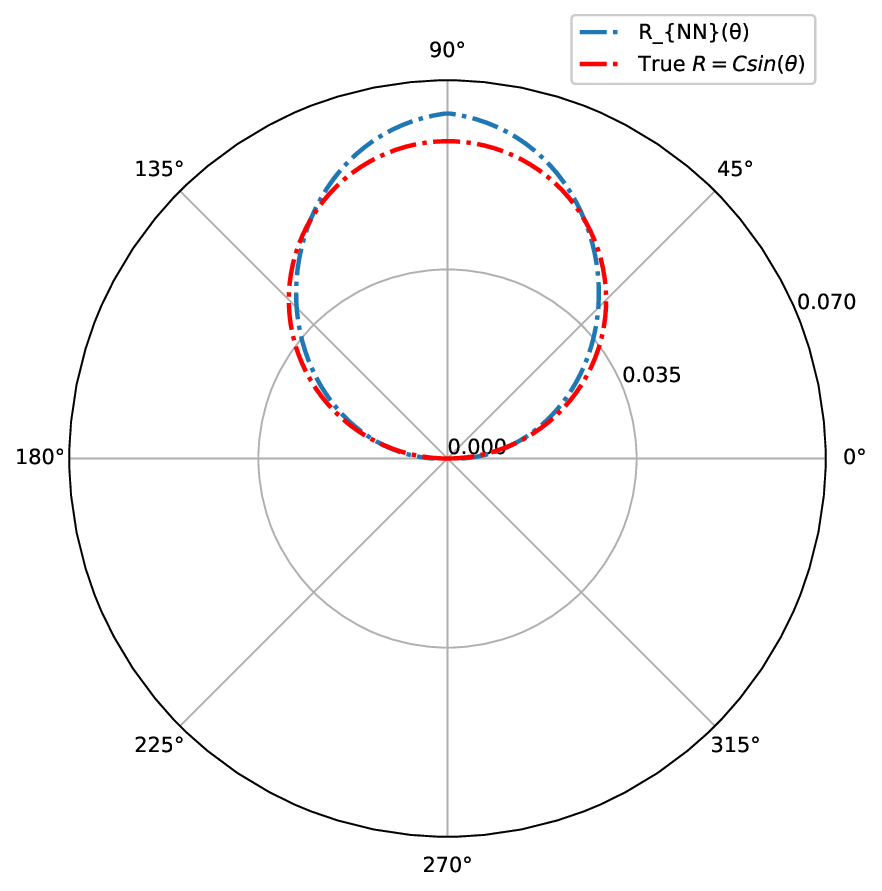}
\label{fig:R-polar}
}
\caption{Learned solution to the stress balance equation using PINN approximation.
The physical parameters used are: surface tension $\si=7.28\times10^{-2}$ Pa$\cdot$m, ambient pressure $p_{\infty,*} = 1.013\times 10^5$ Pa, and bubble volume $V=5\times10^{-4}$ m$^3$.
The penalty weights are $\la_{\rm SB} = 10^3$, $\la_V = 1$, $\la_B = 10^{-6}$, and $\la_S = 10^3$.}
\end{figure}

For training, we use a fully connected feedforward neural network with four layers: three hidden layers of 50 nodes each and a final output layer.
The activation function is {\bf Tanh()}. 
The neural network is trained using the {\bf Adam} optimizer with a learning rate of $10^{-4}$ for $10000$ epochs.
The simulation is implemented in Python using the \texttt{PyTorch} library, which provides automatic differentiation (via \texttt{autograd}) and tools for constructing and training neural network architectures.
Figure~(\ref{fig:3D}) shows the three-dimensional surface of the learned equilibrium bubble, while Figures~(\ref{fig:R}) and (\ref{fig:R-polar}) compare the learned radial profile $R_{\rm NN}(\th)$ against the ground truth $R(\th) = C \sin \th$ in Cartesian and polar coordinates, respectively.
The results demonstrate that the PINN-based approximation captures the overall shape of the exact solution, with a relative root mean square error (rRMSE) of $5.481 \times 10^{-2}$.

\renewcommand{\thesubfigure}{\alph{subfigure}}

 \captionsetup[subfloat]{captionskip=5pt, width=0.3\linewidth, labelformat=parens, labelfont=normalfont}

\FloatBarrier
\FloatBarrier

\bigskip

\end{document}